\documentclass[11pt]{article}

\RequirePackage[OT1]{fontenc}
\RequirePackage{amsthm,amsmath,amssymb,comment}
\usepackage{}
\RequirePackage[numbers]{natbib}
\RequirePackage[colorlinks,citecolor=blue,urlcolor=blue]{hyperref}

\numberwithin{equation}{section}
\theoremstyle{plain}
\newtheorem{theorem}{Theorem}[section]
\newtheorem{definition}{Definition}[section]
\newtheorem{condition}{Condition}[section]
\newtheorem{lemma}{Lemma}[section]
\def\risk{\mathit{R}}
\def\K{\mathcal{K}}
\def\MR{\text{MR}}

\date{}

\title{Ordered Smoothers With Exponential Weighting\thanks{This research was supported  by
Laboratory of Structural Methods of Predictive Modeling and Optimization, Moscow Institute Physics and Technology, Russian Federation government grant, ag. 11.G34.31.0073.}}
\author{
Chernousova, E.,\thanks{
 Moscow Institute of Physics and Technology,
Institutski per. 9, Dolgoprudny, 141700, Russsia, \texttt{lena-ezhova@rambler.ru}
}\quad
Golubev, Yu.\thanks{CNRS and Institute for  Information Transmission Problems,
Bolshoy Karetny per. 19,  Moscow, 127994, 
Russia, \texttt{golubev.yuri@gmail.com}}
\quad and Krymova, E.\thanks{DATADVANCE and  Institute for  Information Transmission Problems, Bolshoy Karetny per. 19, Moscow, 127994, Russia,
\texttt{ekkrym@gmail.com}
 }
}

\begin{document}
\maketitle

\begin{abstract}
The main goal in this paper is to propose a new method for deriving oracle inequalities related to the exponential weighting method. For the sake of simplicity we  focus on recovering an unknown vector from  noisy data with the help of a family of ordered smoothers. The estimators withing this    family
are aggregated using   the exponential weighting  and the aim is 
to  control the risk of the aggregated estimate.  Based on simple probabilistic properties of the unbiased risk estimate,  we derive new oracle inequalities and show that the exponential weighting permits to improve Kneip's oracle inequality \cite{K}.
\end{abstract}




\section{Introduction and main results}

This paper deals with the simplest linear model
\begin{equation}\label{eq1.1}
Y_{i}=\mu _{i}+\sigma \xi _{i} , \quad i=1,2,\ldots,n,
\end{equation}
where $\xi _{i}$ is a standard white Gaussian noise.  For the sake of simplicity  it is assumed that the noise level $\sigma>0$
is known. 

The goal  is to estimate an unknown vector
 $\mu \in \mathbb{R}^n$  based on the data $Y=(Y_1,\ldots,Y_n)^\top$.
In this paper,  $\mu$ is  recovered with the help  of linear
estimates
\begin{equation}\label{eq1.2}
\hat\mu _{i}^{h}(Y)=h_{i}Y_{i},\ h \in \mathcal{H},
\end{equation}
where $\mathcal{H}$ is a  finite  set of  vectors in $\mathbb{R}^n$ which will be described later on.

In what follows, the risk of an estimate $\hat{\mu}(Y)=(\hat{\mu}_1(Y),\ldots,\hat{\mu}_n(Y))^\top$
is measured by
\[
\risk(\hat{\mu },\mu )=\mathbf{E}_{\mu }\Vert \hat{\mu }(Y)
-\mu \Vert ^{2},
\]
where $\mathbf{E}_{\mu }$ is the expectation  with respect to the measure
 $\mathbf{P}_{\mu}$ generated by the observations from  \eqref{eq1.1} and
$\left\| \cdot \right\|, \ \langle \cdot ,\cdot \rangle$
stand for the  norm and the inner product in $\mathbb{R}^n$
\[
\left\| x\right\| ^{2}=\sum_{i=1}^{n }x_{i}^{2},\quad \langle
x,y\rangle =\sum_{i=1}^{n}x_{i}y_{i}.
\]

Since the mean square risk of $\hat{\mu}^h(Y)$ 
\[
R(\hat{\mu}^h,\mu) =\|(1-h)\mu\|^2+\sigma^2\|h\|^2
\]
  depends on $h \in \mathcal{H}$,  one can minimize it choosing properly $h\in \mathcal{H}$.  Very often the minimal risk
\begin{equation*}
r^{\mathcal{H}}(\mu) =\min_{h\in \mathcal{H}}R(\hat{\mu}^h,\mu)
\end{equation*}
is called the oracle risk.  

Obviously, one cannot make use of the oracle estimate
\[
\mu^*(Y)=h^*\cdot Y,
\quad
h^*=\arg\min_{h\in \mathcal{H}}R(\hat{\mu}^h,\mu)
\]
 because it depends on the underlying vector. However, one could try 
 to  construct an estimator $\tilde{\mu}^{\mathcal{H}}(Y)$ based on
 the family of linear  estimates
$
\hat{\mu}^h(Y) , \ h\in \mathcal{H},
$
 with the  risk mimicking  the oracle risk. This idea means that the risk of $\tilde{\mu}^{\mathcal{H}}(Y)$ should be bounded by the so-called oracle inequality
\begin{equation*}
R(\tilde{\mu}^\mathcal{H},\mu)\le r^{\mathcal{H}}(\mu)+\tilde{\Delta}^{\mathcal{H}}(\mu),
\end{equation*}
which holds uniformly in $\mu\in\mathbb{R}^n$. 
Heuristically, this inequality assumes that the remainder term $\tilde{\Delta}^{\mathcal{H}}(\mu)$ is smaller than the oracle risk for all $\mu$. In general,  such an estimator doesn't exist, but  for certain statistical models it possible to construct  an estimator $\tilde{\mu}^{\mathcal{H}}(Y)$  (see, e.g., Theorem \ref{th1.1} below) such that:
\begin{itemize}
\item $\tilde{\Delta}^{\mathcal{H}}(\mu)\le \tilde{C}r^{\mathcal{H}}(\mu)$ for all $\mu\in \mathbb{R}^n$, where $\tilde{C}>1$ is a constant.
\item $\tilde{\Delta}^{\mathcal{H}}(\mu)\ll r^{\mathcal{H}}(\mu)$ for 
all $\mu : \ r^{\mathcal{H}}(\mu)\gg \sigma^2$.
\end{itemize}     
It is  well-known that one can find the estimator with the above properties
 provided that $\mathcal{H}$ is not very rich (see, e.g., \cite{BM}). In particular, as shown in \cite{K}, this can be done for the so-called ordered smoothers.  This is why  this paper deals with
    $\mathcal{H}$ containing solely ordered multipliers defined as follows:
\begin{definition}\label{de1.1} $\mathcal{H}$ is a set of ordered multipliers if
\begin{itemize}
\item  $h_i\in [0,1],\ i=1,\ldots,n $ for all $h\in \mathcal{H}$,
\item   $h_{i+1}\le h_{i}, \ i=1,\ldots,n$ for all $h\in \mathcal{H}$,
\item  if for some integer $k$ and some $h,g\in \mathcal{H}$,
 $h_k < g_k$, then $ h_i\le g_i$ for all $ i=1,\ldots,n$.
\end{itemize}
\end{definition}
The last condition means that vectors in $\mathcal{H}$ may be naturally ordered, since for any $h,g\in \mathcal{H}$ there are only two
possibilities $h_i\le g_i$ or $h_i\ge g_i$ for all $i=1,\ldots,n$. Therefore the estimators from \eqref{eq1.2} are often called ordered smoothers \cite{K}.

Notice  that ordered smoothers are  
common  in  statistics (see, e.g., \cite{K}).
Below we give two basic examples, where   these smoothers appear naturally.

\bigskip

\noindent
\textbf{Smoothing splines.} They are usually used in recovering  smooth regression functions $f(x), \ x\in [0,1],$ given  the noisy observations
\begin{equation}\label{eq1.3}
Z_i=f(x_i)+\sigma \xi_i', \quad  i=1,\ldots,n,\
\end{equation}
where $x_i\in (0,1)$ and $\xi_i'$  are i.i.d. Gaussian random variables with zero mean and unit variance. It is well known that smoothing spline is defined by
\begin{equation}\label{eq1.4}
\hat{f}_\alpha(x,Z)=\arg\min_{f}\biggl\{\sum_{i=1}^n [Z_i-f(x_i)]^2+\alpha \int_0^1 [f^{(m)}(x)]^2\biggr\},
\end{equation}
where $f^{(m)}(\cdot)$ denotes the derivative of order $m$ and  $\alpha>0$ is a smoothing parameter which is usually chosen with the help of the Generalized Cross Validation  (see, e.g., \cite{W}). 

 To transform this model into the sequence space model \eqref{eq1.1},  consider the Demmler-Reinsch \cite{DR} basis $\psi_k(x),\ x\in[0,1], \ k=1,\ldots,n$ having double orthogonality
\begin{equation*}
\begin{split}
\langle \psi_k,\psi_l\rangle_n=\delta_{kl},\quad
\int_0^1 \psi_k^{(m)}(x) \psi_{l}^{(m)}(x)\, dx =\delta_{kl} \lambda_k,\ k,l=1,\ldots,n,
\end{split}
 \end{equation*}
where here and below $\langle u, v \rangle_n$ stands for the
inner product
\[
\langle u,v\rangle_n=\sum_{i=1}^n u(x_i)v(x_i).
\]
It is assumed  for definiteness that the eigenvalues $\lambda_k$ are  sorted in ascending order 
$
\lambda_1\le \ldots\le \lambda_n.
$

With this basis we can represent the underlying function as follows:
\begin{equation}\label{eq1.5}
f(x)=\sum_{k=1}^n \psi_k(x)\mu_k
\end{equation}
and we get from \eqref{eq1.3}
\begin{equation*}
Y_k=\langle Z, \psi_k\rangle_n=\mu_k+\sigma \xi_k.
\end{equation*}

Next,  substituting \eqref{eq1.5} in \eqref{eq1.4},
we arrive at
\begin{equation*}\label{i:eq8}
\hat{f}_\alpha(x,Y) = \sum_{k=1}^n \hat{\mu}_k \psi_{k}(x),
\end{equation*}
where
\[
\hat{\mu}=\arg\min_{\mu}\biggl\{\sum_{k=1}^n[Y_k-\mu_k]^2
+\alpha \sum_{k=1}^n \lambda_k \mu_k^2\biggr\}.
\]
It is seen easily that
\begin{equation*}
\hat{\mu}_k =\frac{Y_k}{1+\alpha\lambda_k}
\end{equation*}
and
thus, we conclude that the model \eqref{eq1.1}--\eqref{eq1.2} is equivalent  to \eqref{eq1.3}--\eqref{eq1.4}  with
\[
h_k=h_k^\alpha=\frac{1}{1+\alpha\lambda_k}.
\]

Notice that a similar  equivalence with 
\[
h_k=\max(1-\alpha \lambda_k,0)
\] takes place in the minimax estimation of smooth regression functions from Sobolev's balls \cite{S}.  

The Demmler-Reinsch basis   is  a very useful tool for statistical analysis of  spline methods.  In practice, this  basis is rarely used since there are very fast algorithms for computing  smoothing  splines (see, e.g., 
\cite{GS} and  \cite{W}).

\bigskip

\noindent
\textbf{Spectral regularizations of large  linear models.} Very often in linear models, we are interested in estimating $X\mu\in\mathbb{R}^n$ based on the observations
\begin{equation}\label{eq1.6}
Z=X\mu +\sigma\xi,
\end{equation}
where $X$ is a known $n\times p$\,-\,matrix and $\xi$ is a standard white Gaussian noise. It is well known that if $X^\top X$ has a large condition number or $p$ is large, then the standard maximum likelihood estimate $X\hat{\mu}^0(Z)$, where
\[
\hat{\mu}^0(Z) =\arg\min_{\mu} \|Z-X\mu\|^2=(X^\top X)^{-1}X^\top Z
\]
may result in a large risk. More precisely, if $X^\top X>0$, then
\[
\mathbf{E}\|X\mu-X\hat{\mu}^0\|^2=\sigma^2p.
\]

Usually the  risk of $X\hat{\mu}^0(Z)$  may be improved with the help of some regularizations. For instance, one can use the Phillps-Tikhonov  regularization \cite{TA}
\begin{equation*}
\hat{\mu}^{\alpha}(Z) =\arg\min_{\mu}\Bigl\{ \|Z-X\mu\|^2+\alpha\|\mu\|^2\Bigr\},
\end{equation*}
where $\alpha>0$ is a smoothing parameter. It is seen easily that
\[
\hat{\mu}^{\alpha}(Z)=[I+\alpha (X^\top X)^{-1}]^{-1}\hat{\mu}^0(Z).
\]

This formula is a particular case of  the so-called spectral regularizations defined as follows (see, e.g., \cite{Engl}):
\begin{equation*}
\hat{\mu}^{\alpha}(Z)=H^\alpha(X^\top X)\hat{\mu}^0(Z),
\end{equation*}
where $H^\alpha(\cdot): \ \mathbb{R}^+\rightarrow [0,1] $ is a 
  function depending on a smoothing parameter $\alpha\in \mathbb{R}^+$. The matrix $H^\alpha(X^\top X)$ may be easily defined 
when $H^\alpha(\lambda),\ \lambda\in \mathbb{R}^+$ admits the Taylor expansion
\[
H^\alpha(\lambda)=\sum_{s=0}^\infty h_s^\alpha\lambda^s.
\]
Then
\[
H^\alpha(X^\top X)=h_0^\alpha I +\sum_{s=1}^\infty h_s^\alpha (X^\top X)^s,
\]
where $I$ is the identity matrix.

Notice that for the Phillps-Tikhonov method  we have
\[
H^\alpha(\lambda)=\frac{1}{1+\alpha/\lambda}, \ \lambda,\, \alpha \in\mathbb{R}^+
\]
and it is clear that this family of functions is ordered in the sense of Definition \ref{de1.1}. Along with the Phillps-Tikhonov regularization, the spectral cut-off and Landweber's iterations (see, e.g., \cite{Engl} for details) are typical examples of ordered smoothers.    

The standard way to construct an equivalent model of the spectral regularizations is to make use of the SVD. Let $e_k, \ k=1,\ldots,p$ and $\lambda_{1}\le \lambda_2\le \ldots\le \lambda_p$ be  eigenvectors and eigenvalues of $X^\top X$.
It is easy to check that
\[
e^*_k=\frac{X e_k}{\sqrt{\lambda_k}}
\]
is an orthonormal basis in $\mathbb{R}^n$.
   Therefore $Z$ defined by  \eqref{eq1.6}  can be represented   in the following equivalent form
\begin{equation}\label{eq1.7}
Y_k^*=\langle e^*_k, Z\rangle=\langle e^*_k,X\mu \rangle +\sigma\xi_k^*,
\end{equation}
where $\xi_k^*$ are i.i.d. $\mathcal{N}(0,1)$. Notice also that
\[
X\hat{\mu}^\alpha(Z)=X H^\alpha(X^\top X)(X^\top X)^{-1} X^\top Z
\]
and hence
\begin{equation}\label{eq1.8}
\begin{split}
\langle X\hat{\mu}^\alpha(Z),e^*_k\rangle=\sum_{s=1}^pY_s \langle X H^\alpha(X^\top X)(X^\top X)^{-1} X^\top e_s^*,e^*_k\rangle\\
=\sum_{s=1}^p Y_s {\lambda_k}\langle H^\alpha(X^\top X)(X^\top X)^{-1}  e_s,e_k\rangle=H^\alpha(\lambda_k)Y_k.
\end{split}
\end{equation}
In view of \eqref{eq1.7} and \eqref{eq1.8}, we see that the spectral regularization methods are equivalent to the statistical model defined by
\eqref{eq1.1} and \eqref{eq1.2}.

\medskip

Nowadays, there are a lot of approaches aimed to construct  estimates  mimicking the oracle risk. At the best of our knowledge, the principal idea in obtaining such estimates goes back to \cite{A} and \cite{M} and related to the method of the unbiased risk estimation \cite{Stein}. The literature on this approach 
is so vast that it would be impractical to cite it here. We mention solely
the following result  by Kneip  \cite{K} since it  plays an  important role in our presentation. Denote by
\begin{equation} \label{eq1.10}
\bar{r}(Y,\hat{\mu}^h)\stackrel{\rm def}{=} \|Y-\hat{\mu}^h(Y)\|^2 +2 \sigma^2\sum_{i=1}^nh_i-\sigma^2n,
\end{equation}
the unbiased risk estimate  of $\hat{\mu}^h(Y)$.
\begin{theorem} \label{th1.1}  Let
\[
\hat{h}=\arg\min_{h\in \mathcal{H}} \bar{r}(Y,\hat{\mu}^h) 
\]
be the minimizer of the unbiased risk estimate.
Then
uniformly in $\mu\in\mathbb{R}^n$,
\begin{equation}\label{eq1.9}
\mathbf{E}_\mu\|\hat{\mu}^{\hat{h}}-\mu \|^2\le r^{\mathcal{H}}(\mu)+K\sigma^2\sqrt{1+\frac{r^{\mathcal{H}}(\mu)}{\sigma^2}},
\end{equation}
where $K$ is a universal constant.
\end{theorem}

Another idea to construct a good estimator based on the family $\hat{\mu}^h,\ h\in \mathcal{H}$ is to aggregate the estimates within this family using a held-out sample. Apparently, this  approach  was firstly developed by Nemirovsky in \cite{N} and independently by Catoni (see \cite{C} for a summary). Later, the method was extended to several statistical models (see, e.g., \cite{Y}, \cite{RT07}, \cite{L}).

To overcome the well-know  drawbacks of sample splitting one would like to aggregate estimators using the same observations for constructing 
estimators and performing the aggregation.  This  can be done, for instance, with the help of the  exponential weighting. The motivation of this method is  related  to the problem of functional aggregation,   see \cite{RT}.
It has been shown that this method yields rather good oracle inequalities for certain statistical models \cite{LB}, \cite{DS}, \cite{RT}.

In context of the considered statistical model,  the exponential weighting  estimate is defined as follows:
\begin{equation*}
\bar{\mu}(Y) =\sum_{h\in \mathcal{H}} w^{h}(Y)\hat{\mu}^h(Y),
\end{equation*}
where 
\begin{equation*}
w^h(Y)=\pi^h\exp\biggl[-\frac{\bar{r}(Y,\hat{\mu}^h)}{2\beta\sigma^2}\biggr]\biggl/
\sum_{g\in \mathcal{H}}\pi^g
\exp\biggl[-\frac{\bar{r}(Y,\hat{\mu}^g)}{2\beta\sigma^2}\biggr], \ \beta>0,
\end{equation*}
and
  $\bar{r}(Y,\hat{\mu}^h)$  is the unbiased risk estimate  of $\hat{\mu}^h(Y)$ defined by \eqref{eq1.10}.

It has been shown in \cite{DS} that for this method  the following oracle inequalities hold.
\begin{theorem}\label{th1.2}
If $\beta \ge 4$, then uniformly in $\mu\in\mathbb{R}^n$
\begin{equation}\label{th1.2}
\begin{split}
R(\bar{\mu},\mu)&\le \min_{\lambda_h\ge 0:\|\lambda\|_1=1 }\biggl\{
\sum_{h\in \mathcal{H}} \lambda_h R(\hat{\mu},\mu)+2\sigma^2\beta \K(\lambda,\pi)\biggr\},\\
R(\bar{\mu},\mu)&\le \min_{h\in\mathcal{H} }\biggl\{
  R(\hat{\mu}^h,\mu)+2\sigma^2\beta \log\frac{1}{\pi_h}
\biggr\},
\end{split}
\end{equation}
where $\K(\cdot,\cdot)$ is the Kullback-Leibler divergence
\[
\K(\lambda,\pi)=\sum_{h\in \mathcal{H}}\lambda_h\log\frac{\lambda_h}{\pi_h}.
\]
\end{theorem}
Notice that for projection methods ($h_k\in \{0,1\}$) this theorem holds for $\beta\ge 2$, see \cite{LB}. 

It is clear that if we want to derive from \eqref{th1.2} an oracle inequality similar to  \eqref{eq1.9}, then we have to chose $\pi_h =(\#\mathcal{H})^{-1}$, where
$\#\mathcal{H}$ denotes the cardinality of $\mathcal{H}$, and thus we arrive at 
\begin{equation*}
R(\bar{\mu},\mu)\le 
  r^{\mathcal{H}}(\mu)+2\sigma^2\beta \log(\#\mathcal{H}).
\end{equation*}   
This oracle inequality is good  only when the cardinality of $\mathcal{H}$ is not very large. If we deal with continuous $\mathcal{H}$ like those related to  splines smoothing with continuous smoothing parameter, this inequality is not good. To some extent, this situation may be improved, see Proposition 2 in  \cite{DS}.
However, looking at the oracle inequality this  proposition, unfortunately, one cannot  say that it is better than \eqref{eq1.9}.   

The main goal is this paper is to show that for the exponential weighting we can get oracle inequalities with smaller remainder terms than that one in Theorem \ref{th1.1},  Equation \eqref{eq1.9}.  

In order to attain this goal and to cover $\mathcal{H}$ with low and very hight cardinalities, we make use of the special prior weights
   defined as follows:
\begin{equation}\label{eq1.12}
\pi^h\stackrel{\rm def}{=}
1-\exp\biggl\{ - \frac{\|h^+\|_1-\|h\|_1}{\beta}\biggr\}.
\end{equation}
Here $$h^+=\min\{g\in \mathcal{H}: g>h\}$$  $\pi^{h_{\max}}=1$, where $h^{\rm max}$ is the maximal multiplier in $\mathcal{H}$, and  $\|\cdot\|_1$ stands for the $l_1$-norm in $\mathbb{R}^n$, i.e.,
$$
\|h\|_1=\sum_{i=1}^n |h_i|.
$$

Along with these weights we will need   also the following condition:
\begin{condition}\label{co1.1}
There exist  constants $K_\circ,\ K^\circ $ such that
\begin{align}\label{cond:H}
&\sum_{i=1}^n (h_i^2-g_i^2) \ge K_\circ  \bigr(\|h\|_1-\|g\|_1\bigl)
\ \text{for all}\ h\ge g\  \text{from}\ \mathcal{H}
,\\
\label{eq1.14}
&\|h^+\|^2\le K^\circ\|h\|^2 \ \text{for all}\ h\in \mathcal{H} .
\end{align}

\end{condition}

The next theorem, yielding an upper bound for the mean square risk of $\bar{\mu}(Y)$, is the main result of this paper.

\begin{theorem}\label{th1.3}
Assume that $\beta\ge 4$ and   Conditions \ref{co1.1}  hold. Then, uniformly in $\mu \in \mathbb{R}^n$,
\begin{equation}\label{eq1.15}
\mathbf{E}_\mu\|\bar{\mu}-\mu\|^2\le r^{\mathcal{H}}(\mu)+2\beta\sigma^2 \log \biggl[\frac{r^{\mathcal{H}}(\mu)}{\sigma^2} +\Psi\biggl(\frac{r^{\mathcal{H}}(\mu)}{\sigma^2}\biggr)\biggr],
\end{equation}
 where  $\Psi(x)>0, \ x\ge 0$ is a nondecreasing  function bounded at $0$ and  such that
$$
{\Psi(x)}\le \frac{Cx}{\log(x)}.
$$
\end{theorem}
We finish this section with a short discussion  concerning this theorem.

\medskip
\noindent
\textbf{Remark 1.} The condition $\beta \ge 4$ may be improved   when the multipliers $h\in \mathcal{H}$  take only two values $0$ and $1$. In this case it is sufficient to assume that $\beta \ge 2$ (see \cite{Go}).  

\medskip
\noindent
\textbf{Remark 2.} In contrast to   Proposition 2 in \cite{DS}, the remainder term in  \eqref{eq1.15} does not depend neither the cardinality of $\mathcal{H}$ nor $n$. It has the same structure as Kneip's oracle inequality in Theorem \ref{th1.1}.

\medskip
\noindent
\textbf{Remark 3.}
Comparing \eqref{eq1.15} with  \eqref{eq1.9}, we see that when
\begin{equation*}
\frac{r^{\mathcal{H}}(\mu)}{\sigma^2}\approx 1,
\end{equation*}
then the remainder terms in \eqref{eq1.9} and \eqref{eq1.15} have the same order,
namely, $C\sigma^2$. However, when
\begin{equation*}
\frac{r^{\mathcal{H}}(\mu)}{\sigma^2}\gg 1,
\end{equation*}
we get
\[
2\beta\sigma^2 \log \biggl[\frac{r^{\mathcal{H}}(\mu)}{\sigma^2}+\Psi\biggl(\frac{r^{\mathcal{H}}(\mu)}{\sigma^2}\biggr)\biggr] \ll K\sigma^2\sqrt{1+\frac{r^{\mathcal{H}}(\mu)}{\sigma^2}},
\]
thus showing that the upper bound for the remainder term in the oracle inequality related to  the exponential weighting  is better than that one in Theorem \ref{th1.1}.

\medskip
\noindent
\textbf{Remark 4.} We carried out numerous simulations to compare numerically the remainder terms in  \eqref{eq1.15} and \eqref{eq1.9}  and to find out what $\beta$ is optimal from a practical viewpoint. Below we summarize what we obtained 
for the smoothing splines.
\begin{itemize}
\item Nearly optimal $\beta$ is close to $1$, but unfortunately,  good oracle inequalities are not available for this case. 
\item There is no big difference between the exponential weighting with $\beta=1$ and the classical unbiased risk estimation. Both methods demonstrate almost similar  statistical performance. However, when ${r^{\mathcal{H}}(\mu)}/{\sigma^2}$ is close to $1$, the exponential weighting works usually  better.
  \item It seems to us that the remainder term in the oracle inequality \eqref{eq1.9} is too
large. We couldn't  see the square-root behavior in the simulations. On the other hand, 
the remainder term in \eqref{eq1.15} seems  adequate to   simulation results.    
\end{itemize}
\section{Proofs}
The approach in  the proof of Theorem \ref{th1.3} is based on a combination of  methods for 
deriving oracle inequalities proposed in \cite{LB} and  \cite{Go}. The cornerstone idea is to make
 use of the following property of the unbiased risk estimate: let
\[
\hat{h}=\arg\min_{h\in\mathcal{H}}\bar{r}(Y,\hat{\mu}^h)
\] 
be the minimizer of the unbiased risk estimate, then for any   sufficiently small $\epsilon <1$,  there exists $\hat{h}^\epsilon \ge \hat{h}$ such that 
 with  the probability 1,
\begin{equation*}
\bar{r}(Y,\hat{\mu}^h)-\bar{r}(Y,\hat{\mu}^{\hat {h}})\ge  2\beta\sigma^2\epsilon\bigl[\|h\|^2-\|\hat{h}\|^2\bigr]-2\beta\sigma^2, 
\end{equation*}
for all $h\ge \hat{h}^\epsilon$. This property means that $w^h(Y)$ are exponentially decreasing for large $h$ and  therefore 
we can obtain  the following entropy bound   (see Lemma \ref{le2.3} in the paper)
\[
\sum_{h\in \mathcal{H}}w^h(Y)\log\frac{\pi^h}{w^h}\le \log\biggl[\sum_{h\le \hat{h}^\epsilon}\pi^h+C\epsilon^{-1}\exp\bigl(C\epsilon^{-1}\bigr)\biggr].
\]
Here and in what follows, $C$ denotes  a generic constant. 

Next, we prove the following upper bound  
\[
\mathbf{E}_\mu\|\hat{h}^\epsilon\|^2 \le \frac{Cr^{\mathcal{H}}(\mu)}{(1-C\epsilon)\sigma^2}+\frac{C}{\epsilon},
\] 
 with the help of Lemma 2 in \cite{G} (see Lemma \ref{le2.5} below). Finally, we combine these facts following the main lines in the  proof of Theorem 5  in \cite{LB}. 

\subsection{Auxiliary facts}
The next lemma collects some useful facts about the prior weights $\pi^h$ defined by \eqref{eq1.12}.
\begin{lemma}  Under Condition \ref{co1.1}, for any $h\in \mathcal{H}$, the following assertions hold:
\begin{itemize}\item
\begin{align}\label{eq2.1}
\sum_{g\ge h}\pi^g \exp\biggl\{ - \frac{\|g\|_1}{\beta}\biggr\}=  \exp\biggl\{-\frac{\|h\|_1}{\beta}\biggr\},
\end{align}
\item  there exists a constant $C_\circ$ such that
\begin{equation}
\sum_{g\le h}\pi^g\le C_\circ \|h \|^2+C_\circ, \label{eq2.2}
\end{equation}
\item there exist  constants $\pi_\circ$ and $\pi^\circ$
such that 
\begin{equation}\label{eq2.3}
\pi_\circ\le \sum_{g: \|h\|^2\le\|g\|^2\le\|h\|^2+1}\pi^g \le \pi^\circ.
\end{equation}
\end{itemize}

\end{lemma}
\textit{Proof.} Denote for brevity
\[
S^h=\sum_{g\ge h}\pi^g \exp\biggl\{ - \frac{\|g\|_1-\|h\|_1}{\beta}\biggr\}.
\]
Then we have
\begin{equation*}
\begin{split}
&S^h-S^{h^+}=\pi^h+\exp\bigg\{-\frac{\|h^+\|_1-\|h\|_1}{\beta}\biggr\}\\ & \times \sum_{g\ge h^+}\pi^g \exp\biggl\{ - \frac{\|g\|_1-\|h^+\|_1}{\beta}\biggr\} -\sum_{g\ge h^+}\pi^g \exp\biggl\{ - \frac{\|g\|_1-\|h^+\|_1}{\beta}\biggr\}\\&=
\pi^h-\biggl\{1-\exp\biggl[-\frac{\|h^+\|_1-\|h\|_1}{\beta}\biggr]\biggr\} S^{h^+}.
\end{split}
\end{equation*}
Therefore in view of the definition of $\pi^h$, it is clear that if $S^{h_{max}}=1$, then
$
S^h=S^{h^+}
$,
thus proving \eqref{eq2.1}.

To prove \eqref{eq2.2}, notice that
\[
\pi^g \le \frac{\|g^+\|_1-\|g\|_1}{\beta}
\]
and hence, by Conditions \eqref{cond:H} and \eqref{eq1.14},
\begin{equation*}
\begin{split}
\sum_{g\le h}\pi^g\le \frac{1}{\beta}\sum_{g\le
h}\bigl[\|g^+\|_1-\|g\|_1\bigr]  = \frac{\|h^+\|_1 -
\|h_{min}\|_1}{\beta}\\ \le \frac{\|h^+\|^2 -
\|h_{min}\|^2}{K_{\circ}\beta} \le \frac{K^\circ\|h\|^2 -
\|h_{min}\|^2}{K_{\circ}\beta}.
\end{split}
\end{equation*}

In order to check \eqref{eq2.3}, consider the following subset in $\mathcal{H}$
\[
\mathcal{G}_h=\bigl\{g: \|h\|^2\le\|g\|^2\le\|h\|^2+1 \bigr\}.
\]
Let $g_h$ be the maximal element in $\mathcal{G}_h$. Then there are two possibilities
\begin{itemize}
\item $\|g_h\|^2\le \|h\|^2+1/2$,
\item $\|g_h\|^2\ge \|h\|^2+1/2$.
\end{itemize}
In the first, case we have
\[
 1/2 \le \|g_h^+\|^2-\|g_h\|^2\le 2\bigl[\|g_h^+\|_1-\|g_h\|_1\bigr]
\]
and therefore by \eqref{eq1.12}
\begin{equation}\label{eq2.4}
\begin{split}
\sum_{g\in\mathcal{G}_h}\pi^g \ge \pi^{g_h}\ge 1-\exp\biggl(-
\frac{\|g_h^+\|_1-\|g_h\|_1}{\beta}\biggr)
=1-\exp[-1/(4\beta)].
\end{split}
\end{equation}
In the case, where $ \|h\|^2+1/2 \le \|g_h\|^2 \le \|h\|^2+1$, we
make use of that by the Taylor expansion, for any $g<g_h$
\begin{equation*}
\begin{split}
\pi^g \ge& \frac{\|g^+\|_1-\|g\|_1}{\beta} \exp\biggl(-\frac{
\|g_h\|_1-\|h\|_1}{\beta}\biggr)\\ \ge&
\frac{\|g^+\|_1-\|g\|_1}{\beta}
\exp\biggl(-\frac{\|g_h\|^2-\|h\|^2}{\beta K_{\circ}}\biggr)\ge
\frac{\|g^+\|_1-\|g\|_1}{\beta} \exp\biggl(-\frac{1 }{\beta
K_{\circ}}\biggr)
\end{split}
\end{equation*}
and thus,
\begin{equation*}
\begin{split}
\sum_{g\in\mathcal{G}_h}\pi^g \ge \frac{\|g_h\|_1-\|h\|_1}{\beta}
\exp\biggl(-\frac{1}{\beta K_{\circ}}\biggr)\ge
\frac{\|g_h\|^2-\|h\|^2}{2\beta} \exp\biggl(-\frac{1}{\beta
K_{\circ}}\biggr)
\\
\ge \frac{1}{4\beta} \exp\biggl(-\frac{1}{\beta K_{\circ}}\biggr).
\end{split}
\end{equation*}
This equation together with \eqref{eq2.4}   guaranties that there exists $\pi_\circ$ such that for all $h\in \mathcal{H}$
\[
\pi_\circ\le \sum_{g\in \mathcal{G}_h}\pi^g.
\]

The proof of the inverse inequality
$
 \sum_{g\in \mathcal{G}_h}\pi^g \le \pi^\circ
$
is quite similar to that one of \eqref{eq2.2}.$\quad\blacksquare$

The following lemma is a cornerstone in the proof of Theorem \ref{th1.3}.
\begin{lemma}\label{le2.2}
For $\beta \ge 4$ the risk of $\bar\mu(Y)$ is bounded from above  as follows:
\begin{equation*}
\mathbf{E}_\mu\|\bar\mu(Y)-\mu\|^2\le \mathbf{E}_\mu{\sum_{h\in \mathcal{H}}
w^h(Y)\bar{r}(Y,\hat{\mu}^h)}.
\end{equation*}
\end{lemma}

\textbf{Proof}. It is based essentially on \cite{LB}. Recall that the unbiased risk estimates for $\bar{\mu}_i(Y)$ and $\bar{\mu}^h_i(Y)$ are computed as follows (see, e.g. \cite{Stein})
\begin{equation}\label{eq2.5}
\begin{split}
&\bar{r}(Y_i,\bar{\mu}_i)=[\bar\mu_i(Y) - Y_i]^2 +
2 \sigma^2 \frac{\partial \bar{\mu}_i(Y)} {\partial Y_i} -
\sigma^2 ,\\
&\bar{r}(Y_i,\hat{\mu}_i^h)=[\hat\mu_i^h(Y) - Y_i]^2 +
2 \sigma^2 h_i -
\sigma^2.
\end{split}
\end{equation}
Since $\sum_{h \in \mathcal{H}}w^h=1$, we have
\begin{equation}\label{eq2.6}
\begin{split}
&[\bar\mu_{i}(Y) - Y_{i}]^2 = \sum_{h\in \mathcal{H}}
w^h(Y)[\bar\mu_{i}(Y) - Y_{i}]^2 \\ &=
 \sum_{h\in \mathcal{H}}
w^h(Y)[\bar\mu_{i}(Y) - \hat{\mu}_i^h(Y) + \hat{\mu}_i^h(Y) -
Y_{i}]^2  \\& =
 \sum_{h\in \mathcal{H}} w^h(Y)[\bar\mu_{i}(Y) -
\hat{\mu}_i^h(Y) ]^2 + \sum_{h\in \mathcal{H}} w^h(Y)[
\hat{\mu}_i^h(Y) - Y_{i}]^2\\&\quad  + 2\sum_{h\in \mathcal{H}}
w^h(Y)[\bar\mu_{i}(Y) - \hat{\mu}_i^h(Y)][\hat{\mu}_i^h(Y) - Y_{i}] \\& =
\sum_{h\in \mathcal{H}} w^h(Y)[\bar\mu_{i}(Y) - \hat{\mu}_i^h(Y)
]^2 + \sum_{h\in \mathcal{H}} w^h(Y)[ \hat{\mu}_i^h(Y) - Y_{i}]^2\\ &\quad +
2\sum_{h\in \mathcal{H}} w^h(Y)[\bar\mu_{i}(Y) -
\hat{\mu}_i^h(Y)][\hat{\mu}_i^h(Y) - \bar\mu_{i}(Y) +
\bar\mu_{i}(Y) -
Y_{i}] \\&=
 -\sum_{h\in \mathcal{H}} w^h(Y)[\bar\mu_{i}(Y) -
\hat{\mu}_i^h(Y) ]^2 + \sum_{h\in \mathcal{H}} w^h(Y)[
\hat{\mu}_i^h(Y) - Y_{i}]^2.
\end{split}
\end{equation}

From the definition of $\bar{\mu}(Y)$ we obviously get
\[
\frac{\partial\bar{\mu}_{i}(Y)} {\partial Y_i}  = \sum_{h\in
\mathcal{H}}w^h(Y)\frac{\partial\hat\mu_{i}^h(Y)} {\partial Y_i}
+ \sum_{h\in \mathcal{H}}\frac{\partial w^h(Y)} {\partial
Y_i}\hat\mu_{i}^h(Y)
\]
and combining this equation with \eqref{eq2.6} (see also \eqref{eq2.5}), we arrive at
\begin{equation}\label{eq2.7}
\begin{split}
&\bar{r}(Y_i,\bar{\mu}_i)=[\bar\mu_{i}(Y) - Y_{i}]^2 +
2 \sigma^2 \frac{\partial \bar{\mu}_{i}(Y)} {\partial Y_i} -
\sigma^2  \\&=
\sum_{h\in \mathcal{H}}  w_h(Y)\biggl\{[ \hat{\mu}_i^h(Y) - Y_{i}]^2
+ 2 \sigma^2\frac{\partial\hat\mu_{i}^h(Y)} {\partial Y_i} -
\sigma^2 \\ &\quad -[\bar\mu_{i}(Y) - \hat{\mu}_i^h(Y) ]^2
+2\sigma^2\frac{\partial \log(w^h(Y))} {\partial
Y_i}\hat\mu_{i}^h(Y)\biggr\}  \\&=
\sum_{h\in \mathcal{H}}
w^h(Y)\bar{r}(Y_{i},\hat{\mu}_{i}^h)+\\&\quad +\sum_{h\in \mathcal{H}}
w^h(Y)\biggl\{-[\bar\mu_{i}(Y) -
\hat{\mu}_i^h(Y) ]^2 +2\sigma^2\frac{\partial \log[w^h(Y)]}
{\partial Y_i}\hat\mu_{i}^h(Y)\biggr\} \\
&=
\sum_{h\in \mathcal{H}}
w^h(Y)\bar{r}(Y_{i},\hat{\mu}_{i}^h)+\sum_{h\in \mathcal{H}}
w^h(Y)\biggl\{-[\bar\mu_{i}(Y) -
\hat{\mu}_i^h(Y) ]^2\\ &\quad +2\sigma^2\frac{\partial \log[w^h(Y)]}
{\partial Y_i}[\hat\mu_{i}^h(Y)-\bar{\mu}_i(Y)]\biggr\}.
\end{split}
\end{equation}
In deriving the above equation it was used that
$
\sum_{h\in \mathcal{H}} w^h(Y) = 1$ and hence
$$
\sum_{h\in \mathcal{H}} \frac{\partial w^h(Y)} {\partial Y_{i}
} =  \sum_{h\in \mathcal{H}} w^h(Y)  \frac{\partial \log
w^h(Y)} {\partial Y_{i} } = 0.
$$

To control the second sum at the right-hand of \eqref{eq2.7}, we make use of the following equation
\begin{equation*}
\begin{split}
\log w^h(Y) =
-\frac{\bar{r}(Y,\hat{\mu}^h)}{2\beta\sigma^2}+\log(\pi^h)
-\log\biggl\{\sum_{g\in \mathcal{H}}\pi^g
\exp\biggl[-\frac{\bar{r}(Y,\hat{\mu}^g)}{2\beta\sigma^2}\biggr]\biggr\}.
 \end{split}
\end{equation*}

Therefore
\begin{equation*}
\begin{split}
&\sum_{h\in \mathcal{H}}  w^h(Y)\frac{\partial \log w^h(Y)}
{\partial Y_{i} }[\hat\mu_{i}^h(Y) - \bar{\mu}_{i}(Y)]\\&\qquad =
-\frac{1}{2\beta\sigma^2} \sum_{h\in \mathcal{H}}
w^h(Y)\frac{\partial \bar{r}(Y,\hat{\mu}^h)} {\partial Y_{i}
}[\hat\mu_{i}^h(Y) - \bar{\mu}_{i}(Y)]
.
\end{split}
\end{equation*}
Substituting in the above equation (see \eqref{eq1.10})
 $$
\frac{\partial \bar{r}(Y_i,\hat{\mu}^h_i)} {\partial
Y_{i} } = 2 (1-h_{i} )^2 Y_{i},
$$
we obtain
\begin{equation}\label{eq2.8}
\begin{split}
&\sum_{h\in \mathcal{H}}  w^h(Y)\frac{\partial \log w^h(Y)}
{\partial Y_{i} }[\hat\mu_{i}^h(Y) - \bar{\mu}_{i}(Y)]\\& \quad =
-\frac{1}{\beta\sigma^2} Y_{i}^2\sum_{h\in \mathcal{H}}
w^h(Y)[h_{i} - 1]^2 [h_{i} - \bar{h}_{i}],
\end{split}
\end{equation}
where $$\bar{h}_{i} = \sum_{h\in \mathcal{H}} w^h(Y)h_{i}.$$

Next noticing that
$$
(1-h_{i})^2 = (1-\bar{h}_{i})^2 +
(\bar{h}_{i}-h_{i})^2 +2(1 - \bar{h}_{i})(\bar{h}_{i} - h_{i}),
$$
we have
\begin{equation*}
\begin{split}
-Y_{i}^2 \sum_{h\in \mathcal{H}} w^h(Y)(h_{i} - 1)^2 (h_{i} -
\bar{h}_{i}) = Y_{i}^2(1-\bar{h}_{i})^2\sum_{h\in \mathcal{H}}
w^h(Y)
(\bar{h}_{i} - h_{i} )  \\ +
Y_{i}^2\sum_{h\in \mathcal{H}} w^h(Y)(\bar{h}_{i}-h_{i})^2
(\bar{h}_{i} - h_{i} + 2 - 2\bar{h}_{i})\\  = 2Y_{i}^2\sum_{h\in
\mathcal{H}} w^h(Y)(\bar{h}_{i}-h_{i})^2 \biggl(1 - \frac{h_{i} +
\bar{h}_{i}} {2}\biggr)\\ \le
2\sum_{h\in \mathcal{H}}
w^h(Y)[\bar{\mu}_{i}(Y)-\hat{\mu}_{i}^h(Y)]^2.
\end{split}
\end{equation*}

Combining this equation with \eqref{eq2.6}--\eqref{eq2.8}, we finish the proof. \quad $\blacksquare$

\bigskip
\noindent

\begin{lemma}\label{le2.3} Suppose $\{q^h\le 1, \ h\in \mathcal{H}\}$
is a nonnegative sequence such that
\begin{itemize}
\item for all $ h\ge \tilde{h}$
\[
q^h \le \exp\biggl\{-\epsilon \biggl[\sum_{s=1}^n (h_i^2-\tilde{h}_i^2)\biggr]-1\biggr\}, \quad \epsilon>0.
\]
\item for some $h^*$  such that $\|h^*\|^2\le \|\tilde h\|^2$
\[
 q^g\ge  q_\circ, \ \text{ for all }\ g\in \mathcal{G}_{h^*}=\Bigl \{g\in \mathcal{H}:\ \|h^*\|^2\le \|g\|^2\le \|h^*\|^2+1\Bigr\}.
\]
\end{itemize}
Let
\[
W^h= \pi ^hq^h\biggl[\sum_{g\in \mathcal{H}}\pi^gq^g\biggr]^{-1}.
\]

Then
\[
H(W^h)\stackrel{\rm def}{=}\sum_{h\in \mathcal{H}} W^h\log\frac{\pi^h}{W^h} \le \log\biggl[\sum_{h\le \tilde{h}}\pi^h+\exp[R(\epsilon)]\biggr],
\]
where
\begin{equation*}
\begin{split}
R(\epsilon)=\log \frac{C[1+\epsilon]}{\epsilon}+\frac{C[1+\epsilon-\epsilon\log(q_\circ)]}{\epsilon q_\circ}.
\end{split}
\end{equation*}
\end{lemma}
\textit{Proof.} Decompose $\mathcal{H}$ onto two subsets
\[
\mathcal{Q}=\{h\ge \tilde{h}\}\cup \mathcal{G}_{h^*}, \quad \mathcal{P}=\mathcal{H}\setminus \mathcal{Q}.
\]
Denote for brevity
\[
P=\sum_{h\in \mathcal{P}}\pi^h q^h , \quad  Q=\sum_{h\in\mathcal{Q}}\pi^h q^h.
\]

By convexity of $\log(x)$
\begin{equation}\label{eq2.9}
\begin{split}
&H(W^h)=\frac{P}{P+Q}\sum_{h\in \mathcal{P}} \frac{\pi^h q^h}{P}\log\frac{(P+Q)/P}{q^h/P}\\
&\quad +\frac{Q}{P+Q}\sum_{h\in \mathcal{Q}} \frac{\pi^hq^h}{Q}\log\frac{(P+Q)/Q}{q^h/Q}\\
&\le \frac{P}{P+Q}\log\frac{P+Q}{P} +\frac{Q}{P+Q}\log\frac{P+Q}{Q}+\frac{P}{P+Q}\log\biggl(\sum_{h\in \mathcal{P}}\pi^h\biggr)\\ &\quad+\frac{1}{P+Q} \biggl[\sum_{h\in \mathcal{Q}}\pi^h q^h\log\frac{1}{ q^h}+Q\log(Q)\biggr].
\end{split}
\end{equation}

Next, notice that $x\log(1/x)$ is an increasing function   when $x\in[0,{\rm e}^{-1}]$. Therefore, using that  $1-\exp(-\epsilon)>(1+\epsilon)^{-1}\epsilon$, we get with Condition \eqref{cond:H} and \eqref{eq2.3}
\begin{equation*}
\begin{split}
&\sum_{h\in \mathcal{Q}}\pi^h q^h\log\frac{1}{ q^h}\le\sum_{h\in \mathcal{G}_{h^*} }\pi^h\log\frac{1}{q_\circ}\\&\quad + \sum_{h\ge \tilde{h}}\pi^h \exp\biggl[-\epsilon\sum_{i=1}^n (h_i^2-\tilde{h}_i^2)-1\biggr]
 \biggl[\epsilon\sum_{i=1}^n (h_i^2-\tilde{h}_i^2)+1\biggr]
\\& \le \pi^\circ\log\frac{1}{q_\circ}+ \frac{1}{\rm e} \sum_{h\ge \tilde{h}}\pi^h \exp\bigl[-\epsilon K_\circ (\|h\|_1-\|\tilde{h}\|_1)\bigr]\\&\quad  \times
 \bigl[K_{\circ}\epsilon (\|h\|_1-\|\tilde{h}\|_1)+1\bigr].
\end{split}
\end{equation*}
To continue this inequality, we make use of (see \eqref{eq1.12}) that
$
\pi^h\le \|h^+\|_1-\|h\|_1
$
and thus,
\begin{equation}\label{eq2.10}
\begin{split}
\sum_{h\in \mathcal{Q}}\pi^h q^h\log\frac{1}{ q^h}
\le
C+\frac{1}{\rm e} \sum_{h\ge \tilde{h}} \exp\bigl[-\epsilon K_\circ (\|h\|_1-\|\tilde{h}\|_1)\bigr]\\
 \bigl[K_{\circ}\epsilon \bigl(\|h\|_1-\|\tilde{h}\|_1\bigr)+1\bigr]\bigl(\|h^+\|_1-\|h\|_1\bigr).
\end{split}
\end{equation}
In order to bound from above the right-hand side at this equation,
consider the set $\{h\in \mathcal{H}: \ h\ge \tilde{h}\}$.  We may assume that  $\{h_k, \ k=1,\ldots \}$ in this set are ordered so that $\|h_{k+1}\|_1\ge \|h_k\|_1$. Denote for brevity
\[
S_{i}=\|h_i\|_1-\|\tilde{h}\|_1.
\]
With these notations we can write
\begin{equation}\label{eq2.11}
\begin{split}
\sum_{h\ge \tilde{h}} \exp\bigl[-\epsilon K_\circ (\|h\|_1-\|\tilde{h}\|_1)\bigr]
 \bigl[2\epsilon \bigl(\|h\|_1-\|\tilde{h}\|_1\bigr)+1\bigr]\bigl(\|h^+\|_1-\|h\|_1\bigr)\\ =
\sum_{i\ge 1} \exp\bigl[-\epsilon K_\circ S_i \bigr]\bigr[\epsilon S_i +1\bigl]
(S_{i+1}-S_i).
\end{split}
\end{equation}
Let us check that
\begin{equation}\label{eq2.12}
\max_{S_k,k\geq 1}\sum_{i\ge 1} \exp\bigl[-\epsilon K_\circ S_i \bigr]
[S_{i+1}-S_i]_+\le \frac{C}{\epsilon K_\circ},
\end{equation}
where $[x]_+=\max(0,x)$.
Solving the equation
\[
\frac{\partial\phantom{S_i}}{\partial S_i}
\sum_{i\ge 1} \exp\bigl[-\epsilon K_\circ S_i \bigr]
[S_{i+1}-S_i]_+=0,
\]
we obtain with a simple algebra
\[
S_{i+1}-S_{i}=\frac{\exp\bigl[\epsilon K_\circ(S_i-S_{i-1})\bigr]-1}{\epsilon  K_\circ}.
\]
Hence
\[
\exp(-\epsilon K_\circ S_i) (S_{i+1}-S_{i})=\frac{\exp[-\epsilon K_\circ S_{i-1}]-\exp[-\epsilon K_\circ S_{i}]}{\epsilon K_\circ}
\]
and summing up these equations, we get \eqref{eq2.12}.

Similar arguments may be applied to prove that
\begin{equation}\label{eq2.13}
\max_{S_k,k\geq 1} \sum_{i\ge 1}S_i \exp\bigl[-\epsilon K_\circ S_i \bigr]
[S_{i+1}-S_i]_+\le \frac{C}{\epsilon K_\circ}.
\end{equation}

With Equations \eqref{eq2.10}--\eqref{eq2.13} we get
\begin{equation}\label{eq2.14}
\sum_{h\in \mathcal{Q}}\pi^h q^h\log\frac{1}{ q^h}
\le \frac{C[1+\epsilon-\epsilon\log(q_\circ)]}{\epsilon}
\end{equation}
and similarly
\begin{equation}\label{eq2.15}
\begin{split}
Q=\sum_{h\in \mathcal{Q}} \pi^h q^h\le
\frac{C[1+\epsilon]}{\epsilon} .
\end{split}
\end{equation}

Therefore
\begin{equation}\label{eq2.16}
\log(Q)\le
\log \frac{C[1+\epsilon]}{\epsilon}.
\end{equation}
Denote for brevity
\[
x=\frac{Q}{P+Q}.
\]
Then with \eqref{eq2.14} and \eqref{eq2.16} we arrive at

\begin{equation}\label{eq2.17}
\begin{split}
H(W^h)\le \max_{x\in [0,1]}\biggl\{-x\log(x)-(1-x)\log(1-x)\\ +(1-x)\log\biggl(\sum_{h\in \mathcal{P}}\pi^h\biggr)
+
xR(\epsilon)\biggr\},
\end{split}
\end{equation}
where
\[
R(\epsilon)=\log \frac{C[1+\epsilon]}{\epsilon}+\frac{C[1+\epsilon-\epsilon \log(q_\circ)]}{\epsilon q_\circ}.
\]
It is seen easily that the minimizer   $x^*$ of the right-hand side at \eqref{eq2.17} is a solution to the following equation
\begin{equation*}\label{en8}
\log{\frac{1-x^*}{x^*}}=\log\biggl(\sum_{h\in \mathcal{P}}\pi^h\biggr)-R(\epsilon)
\end{equation*}
and thus
\[
x^*=\biggl\{1+\biggl(\sum_{h\in \mathcal{P}}\pi^h\biggr)\exp[-R(\epsilon)]\biggr\}^{-1}.
\]
Therefore from (\ref{eq2.17}) we get
\begin{equation*}
\begin{split}
H(W^h)\le&\log\biggl(\sum_{h\in\mathcal{P}}\pi^h\biggr)-\log(1-x^*)\\&-x^*\biggl[\log\frac{x^*}{1-x^*}+\log\biggl(\sum_{h\in\mathcal{P}}\pi^h\biggr)-R(\epsilon)\biggr]\\
=&\log\biggl(\sum_{h\in\mathcal{P}}\pi^h\biggr)-\log(1-x^*) =
\log\biggl[\sum_{h\in\mathcal{P}}\pi^h+{\rm
e}^{R(\epsilon)}\biggr]\\ \le&\log\biggl[\sum_{h<\tilde{h}}\pi^h+{\rm
e}^{R(\epsilon)}\biggr].\quad \blacksquare
\end{split}
\end{equation*}

\begin{lemma}
\label{le2.4} Let $\xi_i$ be i.i.d. $\mathcal{N}(0,1)$ and $\mathcal{G}$
be a set of ordered sequences.
Then for any $\alpha>0$
\begin{align}\label{eq2.18}
&\mathbf{E}\max_{g\in \mathcal{G}}\biggl\{\pm \sum_{i=1}^{n}(g_i^2-2g_i)(\xi
_{i}^{2}-1)-{\alpha}\sum_{i=1}^n g_i^2\biggr\} \leq \frac{C}{\alpha },\\ \label{eq2.19}
&\mathbf{E}\max_{g\in \mathcal{G}}\left\{ \sum_{i=1}^n (1-g_i)^2\xi_i\mu_i-{\alpha }\sum_{i=1}^n (1-g_i)^2\mu_i^2\right\} \leq \frac{C
}{\alpha }.
\end{align}
\end{lemma}
\textit{Proof.} It follows from Lemma 2 in \cite{G}.

\begin{lemma}\label{le2.5}
Let
\begin{equation}\label{eq2.20}
\begin{split}
\hat{h}^\epsilon=\max \Bigl\{h: \bigl[\bar{r}(Y,\hat{\mu}^h)- \bar{r}^{\mathcal{H}}(Y)\bigr] \le 2\beta\epsilon\sigma^2\bigl[\|h\|^2-\|\hat{h}\|^2\bigr]+2\beta\sigma^2\Bigr\},
\end{split}
\end{equation}
where $\epsilon\in (0,(5\beta)^{-1})$  and
$$
\bar{r}^{\mathcal{H}}(Y)=\min_{h\in \mathcal{H}}\bar{r}(Y,\hat{\mu}^h),
 \quad \hat{h}=\arg\min_{h\in\mathcal{H}}\bar{r}(Y,\hat{\mu}^h).
$$
Then

\begin{equation}\label{eq2.21}
\begin{split}
\sigma^2\mathbf{E}_\mu\|\hat{h}^\epsilon\|^2 \le
\frac{1}{1-5\beta\epsilon} {r}^{\mathcal{H}}(\mu)+
\frac{C\sigma^2}{(1-5\beta\epsilon)\epsilon}.
\end{split}
\end{equation}
\end{lemma}

\textit{Proof.} By the definition of $\bar{r}(Y,\hat{\mu}^h)$, see \eqref{eq1.10} and \eqref{eq2.20}, we get
\begin{equation*}
\begin{split}
\hat{h}^\epsilon =& \max\biggl\{h:\ \|(1-h)\mu\|^2+
\sigma^2(1-2\beta\epsilon)\|h\|^2\\& +2\sigma
\sum_{i=1}^\infty(1-h_i)^2\mu_i \xi_i+ \sigma^2
\sum_{i=1}^{\infty}(h_i^2-2h_i)(\xi_i^2-1)   \\ \le&
\|(1-\hat{h})\mu\|^2+ \sigma^2(1-2\beta\epsilon)\|\hat{h}\|^2\\
&+2\sigma \sum_{i=1}^\infty(1-\hat{h}_i)^2\mu_i \xi_i+ \sigma^2
\sum_{i=1}^{\infty}(\hat{h}_i^2-2\hat{h}_i)(\xi_i^2-1)
  +2\beta\sigma^2
\biggr\}.
\end{split}
\end{equation*}
    Let us fix some $\gamma\in(0,1)$. Then we can rewrite the above equation
as follows:
\begin{equation*}
\begin{split}
\hat{h}^\epsilon =\max\biggl\{h:\ (1-\gamma)\|(1-h)\mu\|^2  +\sigma^2(1-2\beta \epsilon-\gamma)\|h\|^2\\
+2\sigma \sum_{i=1}^\infty(1-h_i)^2\mu_i \xi_i+\gamma \|(1-h)\mu\|^2 \\+
\sigma^2 \sum_{i=1}^{\infty}(h_i^2-2h_i)(\xi_i^2-1) +\gamma \sigma^2\|h\|^2
  \\
\le (1+\gamma)\|(1-\hat{h})\mu\|^2
+\sigma^2(1-2\beta\epsilon+\gamma)\|\hat{h}\|^2\\ + 2\sigma
\sum_{i=1}^\infty(1-\hat{h}_i)^2\mu_i \xi_i -
\gamma\|(1-\hat{h})\mu\|^2
\\ + \sigma^2 \sum_{i=1}^{\infty}(\hat{h}_i^2-2\hat{h}_i)(\xi_i^2-1)
-\gamma\sigma^2\|\hat{h}\|^2
+2\beta \sigma^2\biggr\}.
\end{split}
\end{equation*}
Therefore
\begin{equation*}
\begin{split}
\hat{h}^\epsilon\le \hat{h'}^\epsilon =&\max\biggl\{h:\ \sigma^2(1-2\beta \epsilon-\gamma)\|h\|^2
\\ &+\min_{g\in \mathcal{H}}\biggl[2\sigma \sum_{i=1}^\infty(1-g_i)^2\mu_i \xi_i+\gamma \|(1-g)\mu\|^2 \biggr]\\
&+\min_{g\in \mathcal{H}}\biggr[\sigma^2 \sum_{i=1}^{\infty}(g_i^2-2g_i)(\xi_i^2-1) +\gamma \sigma^2\|g\|^2\biggr]  \\
\le& (1-2\beta\epsilon+\gamma)\bigr[\|(1-\hat{h})\mu\|^2 +\sigma^2\|\hat{h}\|^2\bigr]\\&+\max_{g\in \mathcal{H}}\biggl[2\sigma \sum_{i=1}^\infty(1-g_i)^2\mu_i \xi_i- (\gamma - 2 \beta\epsilon) \|(1-g)\mu\|^2 \biggr]\\
&+\min_{g\in \mathcal{H}}\biggr[\sigma^2 \sum_{i=1}^{\infty}(g_i^2-2g_i)(\xi_i^2-1) +\gamma \sigma^2\|g\|^2\biggr]
+2\beta\sigma^2\biggr\}.
\end{split}
\end{equation*}

Next, bounding $max$ and $min$ in this equation with the help of Lemma \ref{le2.4}, we arrive at
\begin{equation*}
\begin{split}
(1-2\beta\epsilon-\gamma)\sigma^2\mathbf{E}_\mu\|\hat{h}^\epsilon\|^2
\le (1-2\beta\epsilon+\gamma)
\mathbf{E}_\mu\Bigl\{\|(1-\hat{h})\mu\|^2
+\sigma^2\|\hat{h}\|^2\Bigr\} \\ +
\frac{C\sigma^2}{\gamma- 2 \beta \epsilon}.
\end{split}
\end{equation*}
Hence, choosing  $\gamma=3\beta\epsilon$, we get 
\begin{equation}\label{eq2.22}
\begin{split}
\sigma^2\mathbf{E}_\mu\|\hat{h}^\epsilon\|^2  \le \frac{1
+\beta\epsilon}{1-5\beta\epsilon}
\mathbf{E}_\mu\Bigl\{\|(1-\hat{h})\mu\|^2
+\sigma^2\|\hat{h}\|^2\Bigr\}\\
+\frac{\sigma^2}{(1-5\beta\epsilon)}\biggl[2\beta + \frac
{C}{\beta\epsilon}\biggr].
\end{split}
\end{equation}

To control the expectation at the right-hand side in \eqref{eq2.22}, notice that for any given $g\in \mathcal{H}$ the following inequality
\begin{equation*}
\sum_{i=1}^{n}[1-\hat{h}_i]^2Y_i^2+2\sigma^2 \sum_{i=1}^n\hat{h}_i\le
\sum_{i=1}^{n}[1-{g}_i]^2Y_i^2+2\sigma^2\sum_{i=1}^n g_i
\end{equation*}
holds.
This yields immediately
\begin{equation*}
\begin{split}
\|(1-\hat{h})\mu\|^2 +\sigma^2\|\hat{h}\|^2+
2\sigma\sum_{i=1}^\infty(1-\hat{h}_i)^2\mu_i \xi_i+
\sigma^2 \sum_{i=1}^{\infty}(\hat{h}_i^2-2\hat{h}_i)(\xi_i^2-1)  \\
\le \|(1-{g})\mu\|^2 +\sigma^2\|{g}\|^2+2\sigma
\sum_{i=1}^\infty(1-{g}_i)^2\mu_i \xi_i+ \sigma^2
\sum_{i=1}^{\infty}({g}_i^2-2{g}_i)(\xi_i^2-1) .
\end{split}
\end{equation*}
So, for any $\gamma\in (0,1)$, we get with this equation and Lemma \ref{le2.4}
\begin{equation*}
\begin{split}
(1-\gamma)\mathbf{E}_\mu\Bigl\{\|(1-\hat{h})\mu\|^2 +\sigma^2\|\hat{h}\|^2 \Bigr\}
\le \|(1-{g})\mu\|^2 +\sigma^2\|{g}\|^2\\
+2\sigma \max_{g}\biggl[-\sum_{i=1}^\infty(1-{g}_i)^2\mu_i
\xi_i-\frac{\gamma}{2\sigma}
\sum_{i=1}^\infty(1-{g}_i)^2\mu_i^2 \biggr]\\
+\sigma^2\max_{g}\biggl[ -
\sum_{i=1}^{\infty}{g}_i^2(\xi_i^2-1)-\frac{\gamma}{2}
\sum_{i=1}^ng_i^2 \biggr]
\\+2\sigma^2\max_{g}\biggl[ -
\sum_{i=1}^{\infty}{g}_i(\xi_i^2-1)-\frac{\gamma}{4}
\sum_{i=1}^ng_i^2 \biggr]\\
\le \|(1-{g})\mu\|^2 +\sigma^2\|{g}\|^2+\frac{C\sigma^2}{\gamma} .
\end{split}
\end{equation*}
Next, minimizing the right-hand side in $g \in \mathcal{H}$, we have
\begin{equation*}
\mathbf{E}_\mu\Bigl\{\|(1-\hat{h})\mu\|^2 +\sigma^2\|\hat{h}\|^2 \Bigr\}
\le \frac{1}{1-\gamma} r^{\mathcal{H}}(\mu)  +\frac{C\sigma^2}{(1-\gamma)\gamma}.
\end{equation*}
Choosing in the above display $\gamma =\beta \epsilon$ and substituting thus obtained inequality in \eqref{eq2.22}, we get \eqref{eq2.21}. \quad
$\blacksquare$

\bigskip

\subsection{Proof of  Theorem \ref{th1.3}}

By the definition of $w^h(Y)$ we have
\begin{equation*}
\begin{split}
\log[w^h(Y)]=-\frac{1}{2\sigma^2 \beta}\bar{r}(Y,\hat{\mu}^h)+\log \pi^h -
\log\biggl\{\sum_{g\in \mathcal{H}}\pi^g\exp\biggl[-\frac{\bar{r}(Y,\hat{\mu}^g)}{2\beta\sigma^2}\biggr]\biggr\}\\
= \frac{1}{2\sigma^2 \beta} \bar{r}(Y,\hat{\mu}^{\hat{h}}) -\frac{1}{2\sigma^2 \beta}\bar{r}(Y,\hat{\mu}^h)+\log \pi^h \\ -
\log\biggl\{\sum_{g\in \mathcal{H}}\pi^g\exp\biggl[-\frac{\bar{r}(Y,\hat{\mu}^g)- \bar{r}(Y,\hat{\mu}^{\hat{h}})}{2\beta\sigma^2}\biggr]\biggr\},
\end{split}
\end{equation*}
where
$
\hat{h}=\arg\min_{h\in \mathcal{H}}\bar{r}(Y,\hat{\mu}^h).
$
Therefore
\begin{equation}\label{eq2.23}
\begin{split}
\sum_{h\in \mathcal{H}}w^h(Y)\bar{r}(Y,\hat{\mu}^h)= \bar{r}(Y,\hat{\mu}^{\hat{h}}) +2\beta\sigma^2 \sum_{h\in \mathcal{H}}w^h(Y)\log\frac{\pi^h}{w^h(Y)}\\
-2\beta\sigma^2\log\biggl\{\sum_{g\in \mathcal{H}}\pi^g\exp\biggl[-\frac{\bar{r}(Y,\hat{\mu}^g)- \bar{r}(Y,\hat{\mu}^{\hat{h}})}{2\beta\sigma^2}\biggr]\biggr\}.
\end{split}
\end{equation}

Next notice that for all $h$ such that $\|\hat{h}\|^2\le \|h\|^2\le \|\hat{h}\|^2+1$ we have by Condition \eqref{cond:H}
\begin{equation}\label{eq2.24}
\begin{split}
\exp\biggl[-\frac{\bar{r}(Y,\hat{\mu}^h)- \bar{r}(Y,\hat{\mu}^{\hat{h}})}{2\beta\sigma^2}\biggr]=\exp\biggl[-\frac{\|(1-h)Y\|^2- (1-\hat{h})Y\|^2}{2\beta\sigma^2}\\
-\frac{1}{\beta}\sum_{i=1}^n (h_i-\hat{h}_i)\biggr]\ge \exp\biggl[-\frac{K_\circ}{\beta}\biggr].
\end{split}
\end{equation}

We begin to control the right-hand side at \eqref{eq2.23} with the last term. Ordering the elements in $\mathcal{H}$, we obtain
\begin{equation}\label{eq2.25}
\begin{split}
&\log\biggl\{\sum_{g\in \mathcal{H}}\pi^g\exp\biggl[-\frac{\bar{r}(Y,\hat{\mu}^g)- \bar{r}(Y,\hat{\mu}^{\hat{h}})}{2\beta\sigma^2}\biggr]\biggr\}\\
&\quad  \ge\log\biggl\{\sum_{g\ge \hat{h}}\pi^g\exp\biggl[-\frac{\bar{r}(Y,\hat{\mu}^g)- \bar{r}(Y,\hat{\mu}^{\hat{h}})}{2\beta\sigma^2}\biggr]\biggr\} \\
&\quad= \log\biggl\{\sum_{g\ge
\hat{h}}\pi^g\exp\biggl[-\frac{\|(1-g)Y\|^2-\|(1-\hat{h})Y\|^2}{2\beta
\sigma^2}\\ &\qquad -
\frac{1}{\beta}\sum_{i=1}^n[g_i-\hat{h}_i]\biggr]\biggr\}\ge
 \log\biggl\{\sum_{g\ge \hat{h}}\pi^g\exp\biggl[ - \frac{1}{\beta}\sum_{i=1}^n[g_i-\hat{h}_i]\biggr]\biggr\}.
 \end{split}
\end{equation}

Thus, from \eqref{eq2.25} and \eqref{eq2.1} we get
\begin{equation}\label{eq2.26}
\log\biggl\{\sum_{g\in \mathcal{H}}\pi^g\exp\biggl[-\frac{\bar{r}(Y,\hat{\mu}^g)- \bar{r}(Y,\hat{\mu}^{\hat{h}})}{2\beta\sigma^2}\biggr]\biggr\}\ge 0.
\end{equation}

Our next step is to bound from above the second term at the right-hand side
of Equation \eqref{eq2.23}.
Lemmas \ref{le2.3} and \ref{le2.5} help us   in solving this problem.

Let $\hat{h}^\epsilon$ be defined by \eqref{eq2.20}. Then for all $h>\hat{h}^\epsilon$
\[
\bigl[\bar{r}(Y,\hat{\mu}^h)- \bar{r}^{\mathcal{H}}(Y)\bigr] \ge
2\beta\epsilon\sigma^2\bigl[\|h\|^2-\|\hat{h}\|^2\bigr]+2\beta\sigma^2
\]
and in view of \eqref{eq2.24} we obtain with Lemma \ref{le2.3} and \eqref{eq2.2}
\begin{equation*}
\begin{split}
\mathbf{E}_\mu\sum_{h\in \mathcal{H}}w^h(Y)\log\frac{\pi^h}{w^h(Y)}\le
\mathbf{E}_\mu\log\biggl\{\sum_{h\le \hat{h}^\epsilon}\pi^h+\exp[R(\epsilon)]\biggr\}\\ \le
\log\Bigl\{\mathbf{E}_\mu\|\hat{h}^\epsilon\|^2+1+\exp[R(\epsilon)]\Bigr\}.
\end{split}
\end{equation*}

Next, we bound   $\mathbf{E}_\mu\|\hat{h}_\epsilon\|^2$ with the help of  \eqref{eq2.21}, thus arriving  at
\begin{equation}\label{eq2.27}
\begin{split}
\mathbf{E}_\mu\sum_{h\in \mathcal{H}} w^h(Y)\bar{r}(Y,\hat{\mu}^h) \le
\mathbf{E}_\mu\bar{r}^{\mathcal{H}}(Y)+ 2\beta\sigma^2 \log\biggl\{
\frac{1}{1-5\beta\epsilon}\times\frac{
r^{\mathcal{H}}(\mu)}{\sigma^2}\\ +
\frac{C}{(1-5\beta\epsilon)\epsilon} +\exp[R(\epsilon)]\biggr\}.
\end{split}
\end{equation}

To finish the proof of the theorem it remains to minimize the right-hand side at this equation in
 $\epsilon$. Assuming that   $ \epsilon\le 1/(5\beta)$ we obtain

\begin{equation*}
\begin{split}
 \frac{1+\beta\epsilon}{1-5\beta\epsilon}\times\frac{ r^{\mathcal{H}}(\mu)}{\sigma^2} +
\frac{C}{(1-5\beta\epsilon)\epsilon}
+\exp[R(\epsilon)]\\
\le \frac{ r^{\mathcal{H}}(\mu)}{\sigma^2}+\frac{ C\epsilon r^{\mathcal{H}}(\mu)}{\sigma^2}+\frac{C}{\epsilon}+\exp\biggl[\frac{C}{\epsilon}\biggr].
\end{split}
\end{equation*}
Therefore choosing
\[
\Psi(x)=C\min_{\epsilon\in[0,1/(5\beta)]}\biggl[\epsilon x +\frac{1}{\epsilon }+\exp\biggl(\frac{C}{\epsilon}\biggr)\biggr]
\]
and combining  Lemma  \ref{le2.2} with \eqref{eq2.27}, we complete the proof of \eqref{eq1.15} since obviously
\[
\mathbf{E}_\mu\bar{r}^{\mathcal{H}}(Y)\le r^{\mathcal{H}}(\mu).
\]

 It is clear that $\Psi(0)$ is bounded from above.
It is also easy to check that as
$\rho\rightarrow 0$,
\begin{equation*}
\begin{split}
\epsilon^*(\rho)=\arg\min_{\epsilon\in[0,1/(5\beta)]}\biggl\{C\epsilon +\rho\biggl[\frac{C}{\epsilon }+\exp\biggl(\frac{C}{\epsilon}\biggr)\biggr]\biggr\}\\
 \approx \log^{-1}\frac{C}{\rho}+2\log^{-2}\frac{C}{\rho}
\log\biggl(\log\frac{C}{\rho}\biggr)
\end{split}
\end{equation*}
and thus,
$
\Psi(x)\le {Cx}/{\log(Cx)}. \quad \blacksquare
$

\bigskip

\end{document}